\documentclass{amsart}
\usepackage[latin1]{inputenc}
\usepackage{amsfonts, amssymb, amsmath, amscd, amsthm}
\usepackage{txfonts, graphicx, color, enumerate}
\usepackage{url}
\usepackage[capitalise]{cleveref}

\theoremstyle{plain}
\newtheorem{thm}{Theorem}[section]

\newtheorem{cor}[thm]{Corollary}

\theoremstyle{definition}
\newtheorem{defn}[thm]{Definition}
\theoremstyle{remark}

\numberwithin{equation}{section}

\begin{document}

\title{Diagonal Knots and the Tau Invariant}

\author{Jackson Arndt}
\author{Malia Jansen}
\author{Payton McBurney}
\author{Katherine Vance}
\address{Simpson College}
\email{katherine.vance@simpson.edu}
\thanks{This research was funded by Dr. Albert H. and Greta A. Bryan through the 2017 Bryan Summer Research Program at Simpson College.  }

\begin{abstract}
In 2003, Ozsv\'ath, Szab\'o, and Rasmussen introduced the $\tau$ invariant for knots, and in 2011, Sarkar published a computational shortcut for the $\tau$ invariant of knots that can be represented by diagonal grid diagrams. Previously, the only knots known to have diagonal grid diagram representations were torus knots.  We prove that all such knots are positive knots, and we produce an example of a knot with a diagonal grid diagram representation which is not a torus knot.
\end{abstract}

\subjclass{57K10; 57K18}
\keywords{knot theory, grid diagrams, tau invariant}

\maketitle

\section{Introduction}

Knot Floer Homology, abbreviated HFK, is a package of knot invariants defined in 2003 by Ozsv\'ath, Szab\'o, and Rasmussen. These invariants can be defined combinatorially using grid diagrams \cite{MR2372850}. The $\tau$ invariant is an integer-valued knot invariant in the HFK package \cite{MR2026543}.  Although computations of $\tau$ are generally very difficult, in 2011 Sarkar gave a computational shortcut for knots represented on diagonal grid diagrams \cite{sarkar2011grid}.  As far as the authors are aware, the only use of this method to date was in the paper in which Sarkar presented it and used it to calculate the $\tau$ invariants of torus knots (as well as gave a new proof of Kronheimer and Mrowka's theorem regarding the unknotting numbers of torus knots).  We ask which knots can be represented using diagonal grid diagrams, prove that for any knot $K$ which can be represented on a diagonal grid diagram, $\tau (K) \geq 0$, and show that not all such knots are torus knots.  

In \cref{sec:prelim} we give background information about grid diagrams and the $\tau$ invariant, as well as describe the computational methods that led us towards our main results.  In \cref{sec:diag_knots} we prove our main results about diagonal knots.

\section{Preliminaries}
\label{sec:prelim}

\begin{figure}
\includegraphics[width=0.75\textwidth]{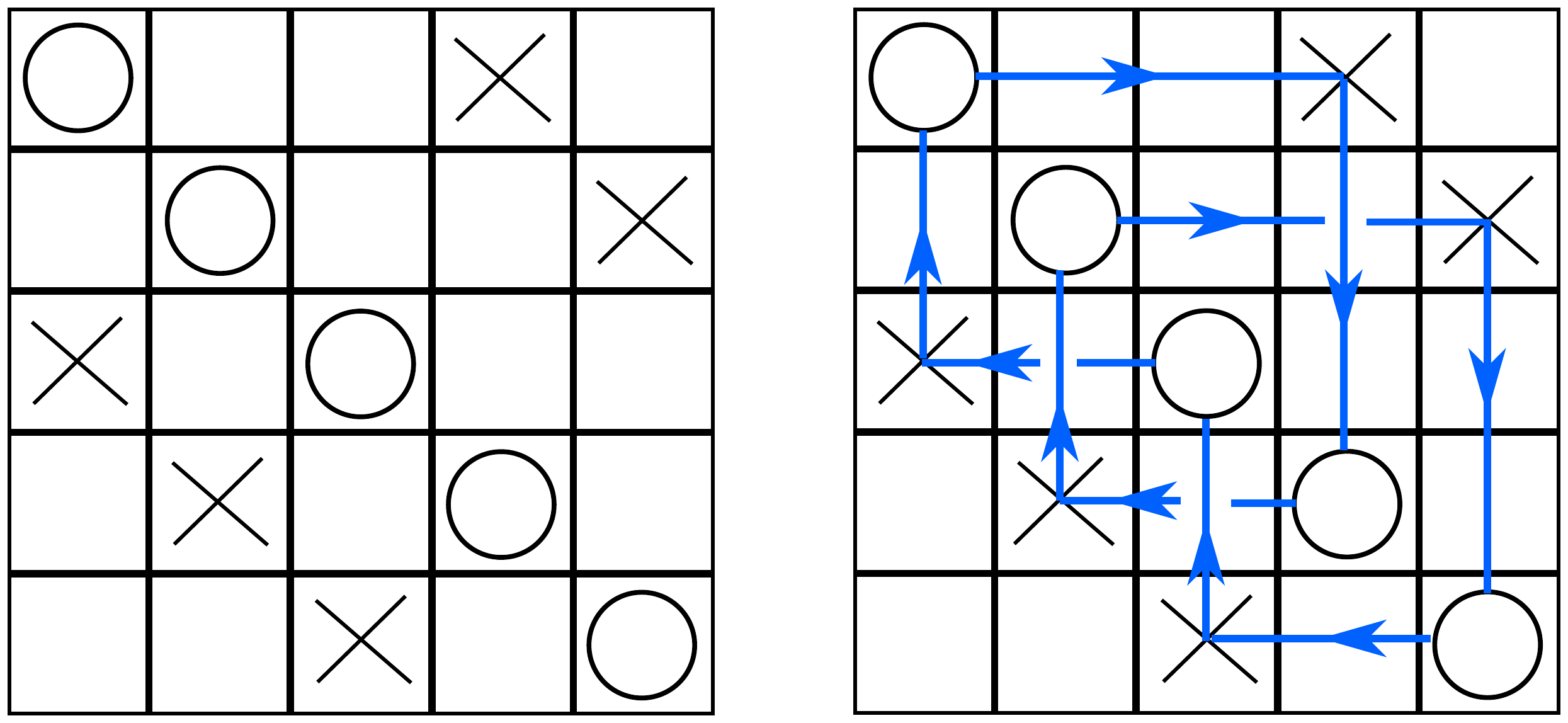}
\caption{A grid diagram representing the trefoil, and the knot diagram recovered.}
\label{grid_trefoil}
\end{figure}

We will use grid diagrams to represent knots.   A grid diagram of size $n$ has $n$ rows and $n$ columns. The rows are numbered bottom to top, while the columns are numbered left to right. In each row and column, there will be exactly one $X$-marking and one $O$-marking, meaning there are $n$ $X$-markings and $n$ $O$-markings. In order to recover the knot given by a diagram, connect the $X$'s to the $O$'s in each column., then connect to the $O$'s to the $X$'s in each row.  At each crossing, the vertical strand passes over the horizontal strand. \cref{grid_trefoil} shows an example of a grid diagram, along with the knot diagram recovered from it.
 
\subsection{The $\tau$ Invariant}
The $\tau$ invariant has some nice geometric properties, including the following bound on slice genus:

\begin{thm}[\cite{MR2026543}] \label{thm:property_tau}
For any knot $K$, $|\tau(K)|\leq g_4(K)$. 
\end{thm}

The $\tau$ invariant is not easy to compute. For some specific cases such as in \cref{thm:sarkar}, there is a straightforward computational method (although for larger grid diagrams, the computation is still quite lengthy). First, we define $\mathcal{J}(A, B)$ as a symmetric, bilinear function on sets in $\mathbb{R}^2$. The function $\mathcal{J}(\{a\}, B)$ is half the number of points in B which are to the top right and bottom left of $a$ (that is, if $a = (a_1,a_2)$ then $\mathcal{J}(\{a\}, B) = \frac{1}{2} |\lbrace (x,y) \in B : (x > a_1 \text{ and } y > a_2) \text{ or } (x < a_1 \text{ and } y < a_2) \rbrace |$), and 
\[
\mathcal{J}(A,B)=\sum_{a\in A} \mathcal{J}(\{a\}, B).
\]

\begin{defn}
A \emph{diagonal knot} is a knot which can be represented on a grid diagram with all of the $O$'s on the diagonal from top left to bottom right. 
\end{defn}

\begin{figure}
\includegraphics[scale=0.3]{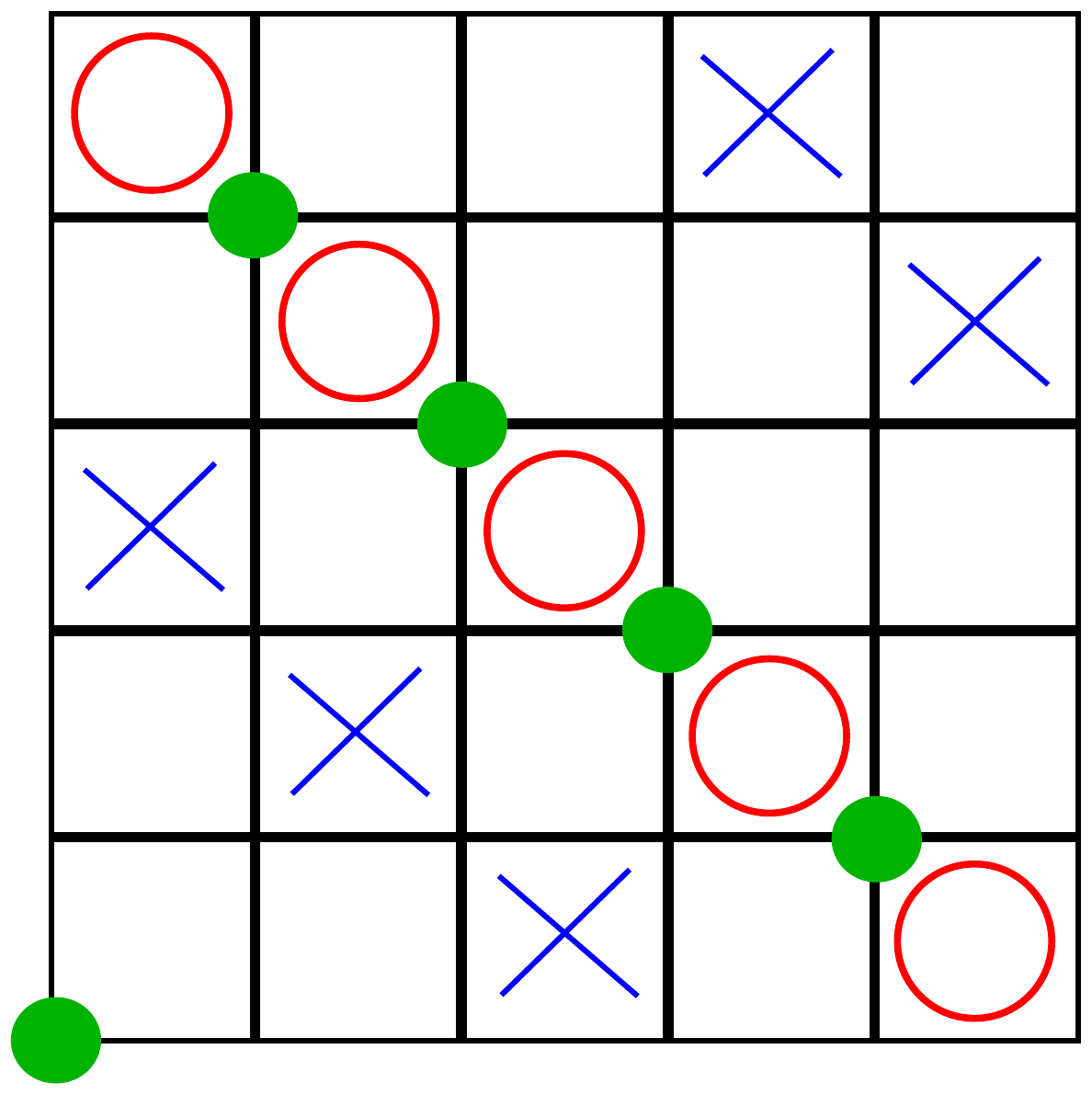}
\caption{An example of a grid diagram where the points used to calculate $\tau$ are indicated with dots.}
\label{color_grid}
\end{figure}

\begin{thm}[\cite{sarkar2011grid}]\label{thm:sarkar} If $K$ is a diagonal knot, then 

\[\tau(K)=\mathcal{J}\left(\mathbf{x}-\frac{1}{2}(\mathbb{X}+\mathbb{O}),(\mathbb{X}-\mathbb{O})\right)-\frac{n-1}{2},
\]

where $\mathbb{X}$ and $\mathbb{O}$ are the sets of $X$- and $O$-markings in a diagonal grid diagram for $K$, $\mathbf{x}$ is the set of points in between the $O$'s on the grid diagram along with one point on the bottom left corner of the grid as shown in \cref{color_grid}, and $n$ is the size of the grid.
\end{thm}

\subsection{Computational methods}

In our investigation of diagonal knot grid diagrams, we wrote code in Java and Mathematica to generate and simplify diagonal knot grid diagrams and to attempt to identify the knots they represent. This code is available on request from the corresponding author (KV). 

\subsubsection{Generating Diagonal Grid Diagrams}
\label{generating_grids}
In order to which knots can be represented using diagonal grid diagrams, we wrote a program that generates all of the X permutations of diagonal knot grid diagrams of a given size, $n$. Once we had generated these we could calculate $\tau$ and other invariants of the knots they represented, with the goal of identifying those knots.

The program generates all X permutations of diagonal grid diagrams which represent knots systematically using recursion. The function takes as its argument the size of the grids diagram we wish to generate, then outputs a tab separated value file containing a list of the X-permutations of all of the diagonal knot grid diagrams of the given size.

\subsubsection{Calculating Tau}
We wrote code in Java that implements Sarkar's method \cite{sarkar2011grid} for calculating the $\tau$ invariant of the knot represented by a given diagonal grid diagram. The program takes just one argument, a permutation specifying the locations of the X-markings in the grid diagram (because this method for computing tau will only work on diagonal grid diagrams, we know that the O's will be on the diagonal).  The program outputs the $\tau$ invariant of the knot represented by the grid diagram with the given X-permutation.  To our knowledge, no other implementation of this algorithm for computing the $\tau$ invariant exists, but one could instead use Ozsv\'ath and Szab\'o's knot Floer homology calculator \cite{HFKcalc}.

\subsubsection{Simplifying Grid Diagrams}
In order to make it easier to identify the knots represented by the diagonal grid diagrams we had generated, we used a program to simplify grid diagrams. Implemented in Java, this program takes as arguments the grid size, X-permutation, and O-permutation of a grid diagram, and returns the X- and O-permutations of the simplified grid diagram.  We started with code written by Park Mikels and Orlando Guerra \cite{parkorlando} that used commutations in order to find adjacent $X$'s and $O$'s in the grid diagram, leading to a destabilization opportunity. Their code, however, did not include cyclic permutations when searching for adjacent $X$'s and $O$'s. We added code that performed cyclic permutations and looped through their code again with the cyclically permuted grid diagram, which slightly increased the number of grid diagrams the program could simplify.  Our code was neither the first nor the best program that performs this task (see, for example, Barbensi and Celoria's GridPyM \cite{MR4722170}), but it was good enough for our purposes.

\subsubsection{Alexander Polynomial}
The last program we wrote calculates the Alexander polynomial of knots given their grid diagram representations. We coded this in Mathematica, using the minesweeper matrix algorithm outlined in \cite{MR2372850} to calculate the Alexander polynomial. The program takes the TSV file outputted but the diagonal knot grid diagram generator described in \cref{generating_grids} and adds a new column containing Alexander polynomials.   Again, our code was neither the first nor the best program to implement this algorithm -- one could instead use SnapPy \cite{SnapPy}, for example.

\section{The class of diagonal knots}
\label{sec:diag_knots}
\subsection{Non-Negativity of $\tau$ for Diagonal Knots}
While using Sarkar's method to calculate $\tau$ for the knots represented by the diagonal grid diagrams generated by our program (described in \cref{generating_grids}), we noticed that all of the computed values are greater than or equal to zero. This observation inspired us to prove \cref{thm:diagpositive}, which implies that the $\tau$ invariant of a diagonal knot is non-negative.  

\begin{figure}
\def\svgwidth{2 in}
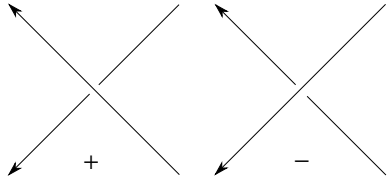
\caption{Examples of positive and negative crossings.}
\label{crossings}
\end{figure}

\begin{defn} \label{positive_crossings}
A crossing in a knot is a \emph{positive crossing} if the crossing follows the right hand rule, which means that if we point our right thumb in the direction of the orientation of the overstrand of the crossing and then sweep our fingers into a fist, our fingers move in the direction of the orientation of the undercrossing strand.  Examples of positive and negative crossings can be seen in \cref{crossings}.
\end{defn} 

\begin{defn} \label{positive_knots}
A knot $K$ is considered to be a \emph{positive knot} if there exists a diagram for $K$ that contains only positive crossings. 
\end{defn}

\begin{figure}
\includegraphics[width=.3\textwidth]{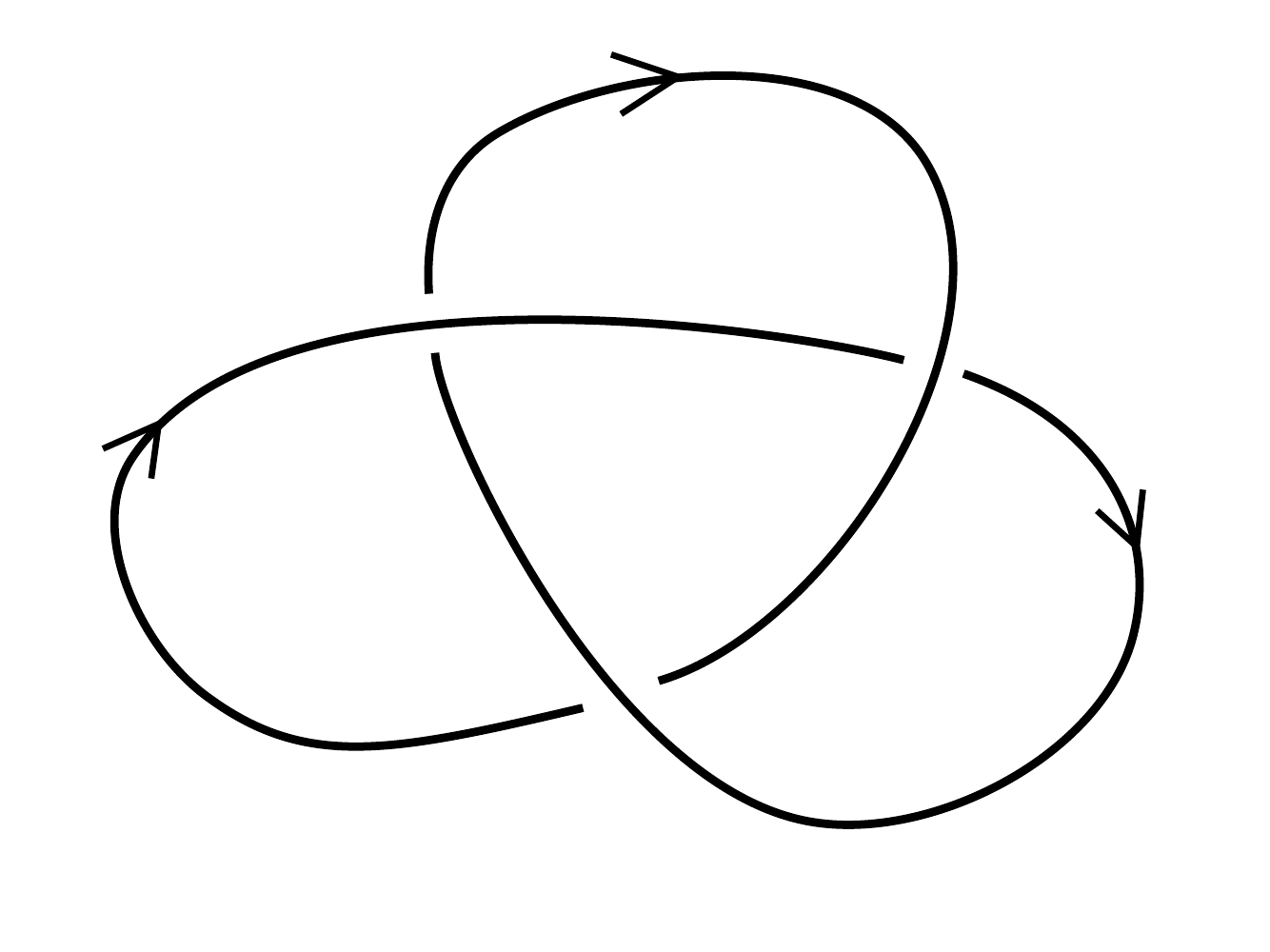}
\caption{The standard diagram for the right handed trefoil contains all positive crossings, and therefore the right-handed trefoil is a positive knot.}
\label{RHT}
\end{figure}

For example, all three of the crossings in the standard diagram of the right handed trefoil  (\cref{RHT}) are positive, so the right-handed trefoil is a positive knot.

\begin{thm}
\label{thm:diagpositive}
Every diagonal knot is a positive knot. 
	\begin{proof}
	A diagonal knot $K$ can be represented using a diagonal grid diagram.  In a grid diagram, the $X$'s connect vertically to the $O$'s and the $O$'s connect horizontally to the $X$'s.  Above the diagonal in the diagonal grid diagram at every crossing the overstrand will be a vertical strand oriented downwards and the understrand will be a horizontal strand oriented from left to right.  Below the diagonal, at every crossing the overstrand will be a vertical strand oriented upwards the the understrand will be a horizontal strand oriented from right to left.  Therefore all of the crossings in a diagonal grid diagram will follow the right hand rule, and so every crossing in the knot diagram for $K$ obtained from the diagonal grid diagram for $K$ is positive.  Thus, $K$ is a positive knot.
	\end{proof}
\end{thm}

\begin{cor}
If $K$ is a diagonal knot, then $\tau(K) \geq 0$.
	\begin{proof}
		By \cref{thm:diagpositive},  $K$ is a positive knot.  Because all positive knots are quasi-positive knots, then $K$ is a quasi-positive knot.  By a result of Livingston, if $K$ is quasi-positive then $\tau(K) = g(K) = g_4(K)$ \cite{livingston2004computations}.  Therefore, $\tau(K) \geq 0$.
	\end{proof}
\end{cor}

\begin{figure} 
\includegraphics[width=1.5 in]{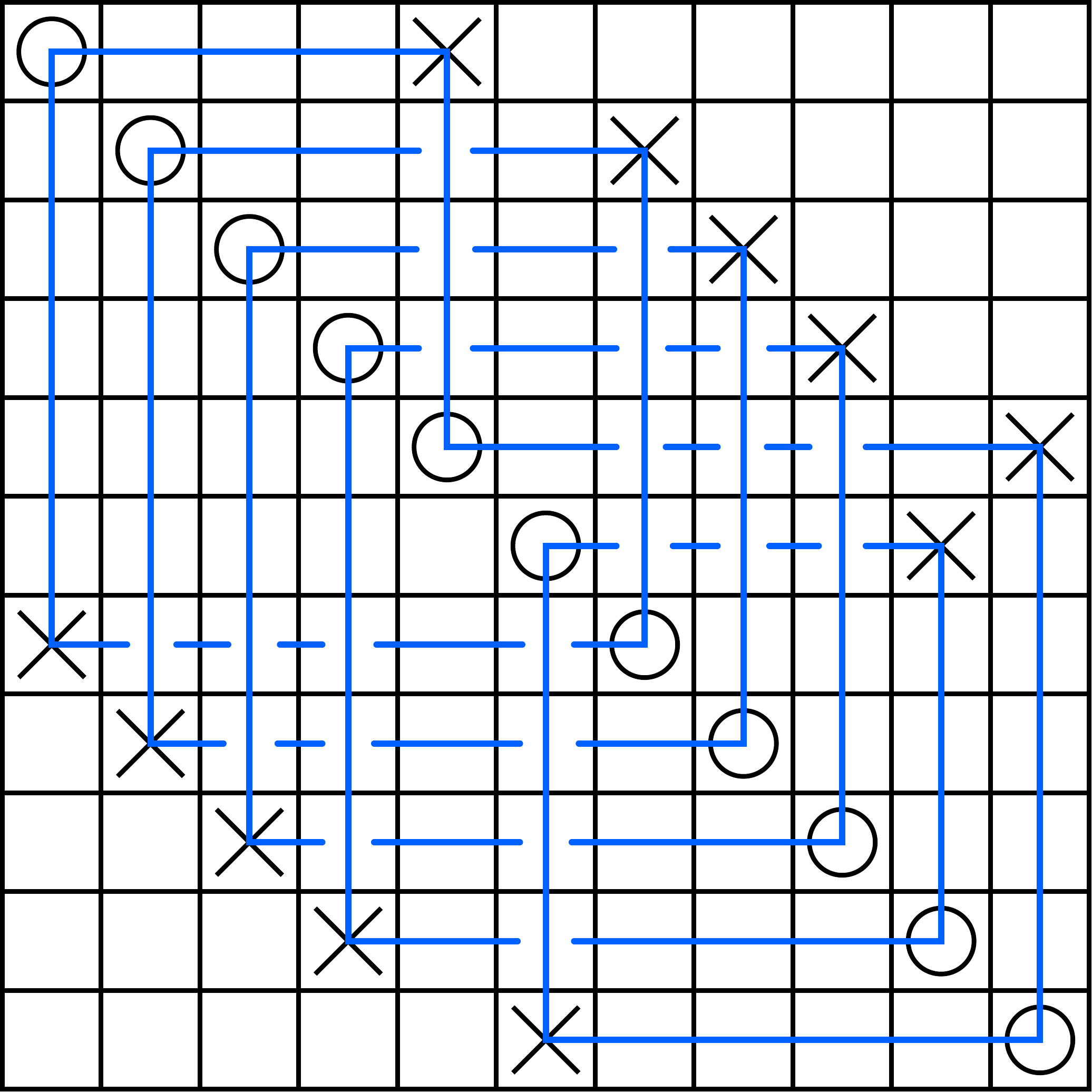}
\caption{A diagonal knot that is not a torus knot.}
\label{twisted_torus_grid}
\end{figure}

\begin{thm}
The set of positive torus knots is a proper subset of the set of diagonal knots.  In particular, there exists a diagonal knot which is neither a torus knot nor a connect sum of torus knots.
	\begin{proof}
	Any positive torus knot $T_{p,q}$ can be represented on a diagonal grid diagram of size $n = p+q$, with $p$ $X$-markings in a diagonal line parallel to and above the $O$-markings and $q$ $X$-markings in a diagonal line parallel to and below the $O$-markings.   \cref{color_grid} shows this diagonal grid diagram representation of the right-handed trefoil, $T_{2,3}$.  Thus, the set of positive torus knots is a subset of the set of diagonal knots.
	
	It remains to show that not all diagonal knots are torus knots.  In the output of the code we wrote to generate diagonal grid diagrams (see \cref{generating_grids}), we discovered a knot represented diagonally with grid size $n = 11$ that is not a torus knot or a connect sum of torus knots. The $X$ permutation for this grid is $\sigma_\mathbb{X} = \{5, 4, 3, 2, 11, 1, 10, 9, 8, 6, 7\}$, and its grid diagram is shown in \cref{twisted_torus_grid}. We used SnapPy to identify this knot as SnapPy census knot m211 \cite{SnapPy}. Through a series of grid diagram commutations, it is straightforward to reduce the number of crossings in the knot diagram from 22 to 21, but it is not obvious whether the number of crossings can be reduced further -- if so, this knot would also appear in Burton's list of knots with up to 19 crossings \cite{zbMATH07760154}. This knot has a $\tau$ invariant of 9 and an Alexander Polynomial of $t^{18} - t^{17} + t^{14} - t^{13} + t^{12} - t^{11} + t^9 - t^7 + t^6 - t^5 + t^4 - t + 1$ (calculated using our Mathematica code and verified with SnapPy \cite{SnapPy}). 

We use two methods to show that the knot represented in \cref{twisted_torus_grid} is not a torus knot.  First, we note that there are eleven torus knots with a grid index less than or equal to eleven, since the grid index of $T_{p,q}$ is $p+q$. Since the grid index of a connect sum of two knots $K$ and $J$ with grid indices $n_K$ and $n_J$, respectively, is $n_K + n_J -1$, there are three connect sums of torus knots with grid index at most eleven (these are $T_{2,3} \# T_{2,3}$, $T_{2,3} \# T_{2,5}$, and $T_{2,3} \# T_{3,4}$).   Comparing the Alexander polynomial of the knot represented in \cref{twisted_torus_grid} to the Alexander polynomials of the torus knots (as listed in \cite{knotatlas}) that can be represented on a size 11 grid diagram, we see that this knot is not a torus knot.  Since the Alexander polynomial of a connect sum of knots is the product of the Alexander polynomials of the summand knots, we compute the Alexander polynomials of $T_{2,3} \# T_{2,3}$, $T_{2,3} \# T_{2,5}$, and $T_{2,3} \# T_{3,4}$ and see that each of them is distinct from the Alexander polynomial of the knot represented in \cref{twisted_torus_grid}.  Therefore this knot is neither a torus knot nor a connect sum of torus knots.  Alternately, we can use SnapPy \cite{SnapPy} to show that this knot is hyperbolic and therefore not a torus knot \cite{MR648524}.  
	\end{proof}
\end{thm}

\section*{Acknowledgments}
This research was funded by Dr. Albert H. and Greta A. Bryan through the 2017 Bryan Summer Research Program at Simpson College.  We are grateful to Park Mikels and Orlando Guerra, who wrote a program as part of their Spring 2017 senior capstone project that was very helpful as we began writing our grid diagram programs.  KV thanks Danmarks Tekniske Universitet (DTU) Compute Institute for Mathematics and Computer Science for welcoming me as a visitor to the Mathematics Section and allowing the use of their library and office resources during the Spring 2024 semester.  We thank the anonymous reviewer for their helpful comments and suggestions.

\bibliographystyle{plain}
\bibliography{KnotBibl}

\end{document}